\documentclass[a4paper, BCOR=0.0mm, DIV=calc, halfparskip, fleqn]{scrartcl}
\usepackage[T1]{fontenc}
\usepackage[utf8]{inputenc}
\usepackage{lmodern}
\usepackage{microtype}
\usepackage[english]{babel}
\usepackage{amsmath, amssymb, amsfonts}
\usepackage{nicefrac}
\newcommand{\diff}[1]{\mathrm{#1}}

\KOMAoptions{twoside=false, twocolumn=false, headinclude=false, footinclude=false, mpinclude=false, pagesize=auto}
\recalctypearea

\begin{document}
\title{On certain definite integrals and infinite series}
\author{Ernst Eduard Kummer, Dr. of Mathematics}
\date{}
\maketitle
\begin{abstract}
We provide a translation of E. E. Kummer's paper "De integralibus quibusdam definitis and seriebus infinitis" \cite{1} \footnote{This paper was translated from Kummer's Latin original "De integralibus quibusdam definitis and seriebus infinitis" by Alexander Aycock}.
\end{abstract}

The definite integrals, that I undertake to treat now, are very closely connected to the infinite series, that I treated in a commentary of this journal  on the hypergeometric series, volume $XV$ page 138 \cite{2} and the following, which, so that they can be represented in a simpler way, I will denote by these functions:\\

\begin{alignat*}{9}
	&1. && 1 && + \frac{\alpha \cdot x}{\beta \cdot 1} && + \frac{\alpha(\alpha+1) \cdot x^2}{\beta(\beta+1) \cdot 1 \cdot 2} && +\frac{\alpha(\alpha+1)(\alpha+2) \cdot x^3}{\beta(\beta+1)(\beta+2) \cdot 1 \cdot 2 \cdot 3} && +\cdots \cdot  && =\varphi{(\alpha, \beta, x)}, &&  &&  \\
	&2. && 1 && + \frac{x}{\alpha \cdot 1} && + \frac{x^2}{\alpha(\alpha +1) \cdot 1 \cdot 2} && + \frac{x^3}{\alpha(\alpha +1)(\alpha +2) \cdot 1 \cdot 2 \cdot 3} && + \cdots \cdot && =\psi{(\alpha,x)}, &&  &&  \\	
	&3.~~ && 1 && - \frac{\alpha \cdot \beta}{1 \cdot x} && + \frac{\alpha(\alpha+1)  \beta(\beta+1)}{1 \cdot 2 \cdot x^2} && - \frac{\alpha(\alpha+1)(\alpha+2)  \beta(\beta+1)(\beta+2)}{1 \cdot 2 \cdot 3 \cdot x^3}&& +\cdots \cdot && = \chi{(\alpha,\beta,x)}. &&  &&   \\	
\end{alignat*}
From this the transformations, found at the cited place, of the series can also be exhibited in this way:\\

\begin{alignat*}{9}
	&4.~~~~~~ && \varphi{(\alpha, \beta, x)} && = e^{x} \cdot \varphi{(\beta-\alpha, \beta, -x)},&&  &&   &&  &&  &&  \\
	&5. && \psi{(\alpha,x)} && = e^{\pm2\sqrt{x}}\varphi{(\alpha-\frac{1}{2},2\alpha-1, \pm4\sqrt{x})}, &&   &&  &&  &&  \\	
\end{alignat*}
which is the same formula as
\[
~~~~~~~~~~~~~~~~6.~~~ \varphi{(\alpha,2\alpha,x)}= e^{\frac{x}{2}}\psi({\alpha+\frac{1}{2},\frac{x^2}{16}})
\]
and
\[
7. ~~~~ \chi{(\alpha, \beta, x)}= \frac{x^{\alpha} \Pi(\beta-\alpha-1)}{\Pi(\beta-1)}\varphi{(\alpha, \alpha-\beta+1,x)}+\frac{x^{\beta} \Pi(\alpha-\beta-1)}{\Pi(\alpha-1)}\varphi{(\beta, \beta-\alpha+1,x)}
\]
After having prepared these things, I want to settle the question about this integral at first
\[
~~~~~~~~~~~~~~~~8. ~~~ y= \int_0^{\infty} u^{\alpha-1} \cdot e^{-u} \cdot e^{-\frac{u}{x}} \diff{d}u,
\]

from where it follows
\[
\frac{\diff{d}y}{\diff{d}x}= -\int_0^{\infty} u^{\alpha-2} \cdot e^{-u} \cdot e^{-\frac{u}{x}} \diff{d}u, ~~\frac{\diff{d}^2y}{\diff{d}x^2}=\int_0^{\infty} u^{\alpha-3} \cdot e^{-u} \cdot e^{-\frac{u}{x}} \diff{d}u,
\]
it is by differentiating of the quantity $u^{\alpha-1} \cdot e^{-u} \cdot e^{-\frac{u}{x}}$:
\[
~~~~~~~~~~~~~~~~~~~~~~~~~\diff{d}(u^{\alpha-1} \cdot e^{-u} \cdot e^{-\frac{u}{x}})
\]
\[
= - u^{\alpha-1} \cdot e^{-u} \cdot e^{-\frac{u}{x}} \diff{d}u+(\alpha-1) u^{\alpha-2} \cdot e^{-u} \cdot e^{-\frac{u}{x}} \diff{d}u+x \cdot  u^{\alpha-3} \cdot e^{-u} \cdot e^{-\frac{u}{x}} \diff{d}u,
\]
and by integration between the boundaries $0$ and $\infty$:

\[
0= -\int_0^{\infty} u^{\alpha-1} \cdot e^{-u} \cdot e^{-\frac{u}{x}} \diff{d}u+(\alpha-1)\int_0^{\infty} u^{\alpha-2} \cdot e^{-u} \cdot e^{-\frac{u}{x}} \diff{d}u
\]
\[
~~~~~~~~~~~~~~~~~~~~~~~~~+x \int_0^{\infty} u^{\alpha-3} \cdot e^{-u} \cdot e^{-\frac{u}{x}} \diff{d}u,
\]
or, what is the same,
\[
~~~~~~~~~~~~~~~~~~~9. ~~~ 0= y+(\alpha-1)\frac{\diff{d}y}{\diff{d}x}-x\frac{\diff{d}^2y}{\diff{d}x^2},
\]
The complete integral of this equation is easily found by means of series, that we denoted by the function $\psi$,
\[
~~~~~~~~~~~~~~~~10. ~~~~ y= A \cdot \psi(1-\alpha, x)+B \cdot x^{\alpha} \psi(1+\alpha,x),
\]
where $A$ and $B$ are arbitrary constants. From there this expression for the presented integral follows
\[
~~~~~\int_0^{\infty} u^{\alpha-1} \cdot e^{-u} \cdot e^{-\frac{u}{x}} \diff{d}u = A \cdot \psi(1-\alpha, x)+B \cdot x^{\alpha} \psi(1+\alpha,x).
\]
The determination of the constant $A$ is easy; for, if we suppose the quantity $\alpha$ to be positive, and put $x=0$, we have
\[
~~~~~~~~~~~~~~~~~~~~~~~\int_0^{\infty} u^{\alpha-1} \cdot e^{-u}  \diff{d}u=A
\]
or
\[
~~~~~~~~~~~~~~~~~~~~~~~~~~~~~~~A=\Pi(\alpha-1).
\]
To determine the constant $B$ in the same way, the integral $y$ has to be transformed by the substitution $u=\frac{x}{v}$, whence it is 
\[
~~~~~~~~\int_0^{\infty} u^{\alpha-1} \cdot e^{-u} \cdot e^{-\frac{u}{x}} \diff{d}u =x^{\alpha}\int_0^{\infty} v^{-\alpha-1} \cdot e^{-v} \cdot e^{-\frac{v}{x}} \diff{d}v, 
\]
after having used this integral transformation, equation (11.) is converted into this one:
\[
~~~\int_0^{\infty} v^{-\alpha-1} \cdot e^{-v} \cdot e^{-\frac{v}{x}} \diff{d}v= A \cdot x^{-\alpha} \psi(1-\alpha,x)+B \cdot \psi(1+\alpha,x),
\]
hence, if we suppose the quantity $\alpha$ to be negative and put $x=0$, we have
\[
~~~~~~~~~~~~~~~~~~~~~~~\int_0^{\infty} v^{-\alpha-1} \cdot e^{-v}  \diff{d}v=B
\]
or
\[
~~~~~~~~~~~~~~~~~~~~~~~~~~~~~~~A=\Pi(-\alpha-1),
\]
after having finally substituted which values of the constants, it is:
\[
~~~~~~~~12. ~~~ \Pi(\alpha-1) \psi(1-\alpha, x)+\Pi(-\alpha-1) x^{\alpha} \psi(1+\alpha,x).
\]
From this determination of the constants certain doubts are to be removed, which can arise from the fact, that the one constant was found, after having put $\alpha>0$, the determination of the other constant on the other hand requires the opposite assumption. But it is nevertheless clear, that these conditions would have been superfluous, if, while determining the constants, we had not used the value $x=0$, but any other positive values, and the values of the constants would not have been other ones. Moreover it is to be beared in mind, that formula (12.) is only valid, if $x$ is a positive quantity, otherwise the integral would become infinite; but if $x$ is positive, this integral has a finite value, whatever the quantity $\alpha$ is, positive or negative.\\

From this formula (12.) one can deduce another integral, which is expressed by two series of the kind $\varphi(\alpha, \beta, x)$. By putting $xv$ in the place of $x$, by multiplying by $e^{-v} \cdot v^{\beta-1} \cdot \diff{d}v$ and integrating between the boundaries $0$ and $\infty$, it is
\[
\int_0^{\infty} \int_0^{\infty} u^{\alpha-1} \cdot e^{-u} \cdot v^{\beta-1} \cdot e^{-v} \cdot e^{-\frac{xv}{u}} \diff{d}u \diff{d}v= \Pi(\alpha-1) \int_0^{\infty} v^{\beta-1} \cdot e^{-v} \cdot \varphi(1-\alpha,xv) \diff{d}v
\]
\[
~~~~~~~~~~~~~~~+\Pi(-\alpha-1)x^{\alpha}\int_0^{\infty} v^{\alpha+\beta-1} \cdot e^{-v} \cdot \varphi(1+\alpha,xv) \diff{d}v,
\]
the integrations with respect to the variable $v$ are easily executed; or it is
\begin{alignat*}{9}
	&&& \int_0^{\infty} v^{\beta-1} e^{-v} \psi(1-\alpha,xv)\diff{d}v && = \Pi(\beta-1)\varphi(\beta,1-\alpha,x),&&  &&   &&  &&  &&  \\
	& && \int_0^{\infty} v^{\alpha+\beta-1} e^{-v} \psi(1+\alpha,xv)\diff{d}v && = \Pi(\alpha+\beta-1)\varphi{(\alpha+\beta,1+\alpha, x)}, &&   &&  &&  &&  \\	
	& && \int_0^{\infty} v^{\beta-1} e^{-v} \cdot e^{-\frac{xv}{u}}\diff{d}v && = \frac{\Pi(\beta-1)}{(1+\frac{x}{u})^{\beta}}, &&   &&  &&  &&  \\
\end{alignat*}
whence
\[
\int_0^{\infty} \int_0^{\infty} u^{\alpha-1} \cdot e^{-u} \cdot v^{\beta-1} \cdot e^{-v} \cdot e^{-\frac{xv}{u}} \diff{d}u \diff{d}v= \Pi(\beta-1)\int_0^{\infty} \frac{u^{\alpha-1} \cdot e^{-u} \diff{d}u}{(1+\frac{x}{u})^{\beta}},
\]
which integral, by putting $ux$ in the place of $u$, is changed into
\[
~~~~~~~~~~~~~~~~~~~~~~~~\Pi(\beta-1)x^{\alpha} \int_0^{\infty} \frac{u^{\alpha+\beta-1}\cdot e^{-ux} \diff{d}u}{(1+u)^{\beta}}
\]
after having substituted which values, we finally have
\[
~~~~~~~~~~~~~~~~~~~~~~~~\Pi(\beta-1)x^{\alpha} \int_0^{\infty}\frac{u^{\alpha+\beta-1} \cdot e^{-ux} \diff{d}u}{(1+u)^{\beta}}
\]
\[
=\Pi(\alpha-1)\Pi(\beta-1)\varphi(\beta,1-\alpha,x)+\Pi(-\alpha-1)\Pi(\alpha+\beta-1)x^{\alpha} \varphi(\alpha+\beta, 1+\alpha,x),
\]
which formula, by changing $\alpha$ into $\alpha-\beta$, obtains this more convenient form
\[
~~~~~~~~~~~~~~~~~~13. ~~~~ \Pi(\beta-1)x^{\alpha} \int_0^{\infty} \frac{u^{\alpha+\beta-1}\cdot e^{-ux} \diff{d}u}{(1+u)^{\beta}}
\]
\[
\frac{\Pi(\alpha-\beta-1)}{\Pi(\alpha-1)}x^{\beta} \cdot \varphi(\beta, \beta-\alpha+1,x)+\frac{\Pi(\beta-\alpha-1)}{\Pi(\beta-1)}x^{\alpha} \cdot \varphi(\alpha, \alpha-\beta+1,x).
\]
Because the one part of this equation, after having interchanged the quantities $\alpha$ and $\beta$, remains the same, it has to be
\[
14. ~~~~\frac{x^{\alpha}}{\Pi(\alpha-1)} \int_0^{\infty} \frac{u^{\alpha-1}\cdot e^{-ux}\cdot \diff{d}u}{(1+u)^{\beta}}=\frac{x^{\beta}}{\Pi(\beta-1)} \int_0^{\infty} \frac{u^{\beta-1}\cdot e^{-ux} \cdot \diff{d}u}{(1+u)^{\alpha}}.
\]
If the transformation, that equation (7.) contains, is applied to formula (13.), it is 
\[
~~~~~~~~~~~~~~~15. ~~~\frac{x^{\alpha}}{\Pi(\alpha-1)} \int_0^{\infty} \frac{u^{\alpha-1}\cdot e^{-ux}\cdot \diff{d}u}{(1+u)^{\beta}}=\chi(\alpha, \beta, x).
\]
Because the series $\chi(\alpha, \beta, x)$ belongs to the class of semiconvergent series, it seems to be neccessary, that formula (15.) is reeinforced by a proof, from which it becomes clear at the same time, that by computation of a certain number of the first terms of this series, the proximate value of this integral is found. For this purpose I use the known equation
\[
1-\frac{\beta}{1}z+\frac{\beta(\beta+1)}{1 \cdot 2}z^2-\cdots \cdot (-1)^{k-1}\frac{\beta(\beta+1)\cdots \cdot (\beta+k-2)}{1 \cdot 2 \cdots \cdot (k-1)}z^{k-1}
\] 
\[
=\frac{1}{(1+z)^{\beta}}-\frac{(-1)^{k}\beta(\beta+1)\cdots \cdot (\beta+k-1)}{1 \cdot 2 \cdots \cdot (k-1)}z^{k} \int_0^1 \frac{(1-u)^k \diff{d}u}{(1+zu)^{\beta+k}},
\]
it is by putting $z=\frac{v}{x}$, by multiplaying by $v^{\alpha-1} \cdot e^{-v} \cdot \diff{d}v$, then by integrating from $v=0$ to $v=\infty$ and dividing by $\Pi(\alpha-1)$
\[
16. ~~~ 1-\dfrac{\alpha \cdot \beta}{1 \cdot x}+\dfrac{\alpha(\alpha+1)\beta(\beta+1)}{1 \cdot 2 \cdot x^2}- \cdots \cdot (-1)^{k-1}\dfrac{\alpha(\alpha+1) \cdots \cdot (\alpha+k-2)\beta(\beta+1) \cdots \cdot (\beta+k-2)}{1 \cdot 2 \cdot 3 \cdots \cdot (k-1) \cdot x^{k-1}}
\]
\[
=\frac{1}{\Pi(\alpha-1)}\int_0^{\infty} \frac{v^{\alpha-1} \cdot e^{-v} \cdot \diff{d}v}{(1+\frac{v}{x})^{\beta}}-\tfrac{(-1)^k\beta(\beta+1) \cdots \cdot (\beta+k-1)}{\Pi(\alpha-1)1 \cdot 2 \cdot 3 \cdots \cdot (k-1) \cdot x^{k}} \int_0^1 \int_0^{\infty} \frac{(1-u)^{k-1} \cdot v^{\alpha+k-1} \cdot e^{-v} \cdot \diff{d}v \cdot \diff{d}u}{(1+\frac{uv}{x})^{\beta+k}},
\]
this double integral along with its coefficients indicates the error, that is committed, if the integral 
\[
\frac{1}{\Pi(\alpha-1)}\int_0^{\infty} \frac{v^{\alpha-1} \cdot e^{-v} \cdot \diff{d}v}{(1+\frac{v}{x})^{\beta}}, ~~~~ \text{or, what is the same} ~~~~ \frac{x^{\alpha}}{\Pi(\alpha-1)}\int_0^{\infty} \frac{v^{\alpha-1} \cdot e^{-vx} \cdot \diff{d}v}{(1+v)^{\beta}}
\]
is computed by the first terms of that series, whose number is $k$, if $k$ is so large, that $\beta+k$ is positive, that quantity, we called the error, changes the sign at the same time as $k$ is converted into $k+1$, or, if a certain number of terms of that series is computed, this sum is either larger or smaller than the desired integral, but if the subsequent term of the series is added, this new sum is smaller than the desired integral, if the sum was larger, if that sum was smaller. Therefore the sums, which that series gives, are alternately too large and too small, and it becomes clear, that the proximate value is found, if the computation is extended to the smallest terms of the semiconvergent series. The same thing can be demonstrated from equation (16.) in this way.\\
Of course it is for positive $\beta+k$:
\[
\int_0^1 \int_0^{\infty} \frac{(1-u)^{k-1} \cdot v^{\alpha+k-1} \cdot e^{-v} \cdot \diff{d}v \diff{d}u}{(1+\frac{uv}{x})^{\beta+k}}<\int_0^1 \int_0^{\infty} (1-u)^{k-1} \cdot v^{\alpha+k-1} \cdot e^{-v} \cdot \diff{d}v \diff{d}u
\]
and
\[
~~~~~~\int_0^1 \int_0^{\infty}(1-u)^{k-1} \cdot v^{\alpha+k-1} \cdot e^{-v} \cdot \diff{d}v \diff{d}u= \frac{\Pi(\alpha+k-1)}{k},
\]

so the error, which is expressed by that double integral, is always smaller than
\[
~~~~~~~~~~~~~~~~~~~~\frac{\beta(\beta+1) \cdots \cdot (\beta+k-1)\Pi(\alpha+k-1)}{1 \cdot 2 \cdot 3 \cdots \cdot (k-1) \cdot \Pi(\alpha-1) x^{k}},
\]
because which is the first term neglected, it follows, that the error is always smaller than that term of the series, to which the summation is extended.\\[2mm]

After having put $\beta=1- \alpha$, equation (15.) is converted into this one:
\[
~~~~~~~~~~~~~~~~~~~~~~~~~\frac{x^{\alpha}}{\Pi(\alpha-1)} \int_0^{\infty} (u+u^2)^{\alpha-1} \cdot e^{-ux} \cdot \diff{d}u
\] 
\[
~~~~~~~=\frac{\Pi(2\alpha-2)}{\Pi(\alpha-1)}x^{1-\alpha} \cdot e^{\frac{x}{2}}\cdot \varphi(1-\alpha,2-2\alpha,x)+ \frac{\Pi(-2\alpha)}{\Pi(-\alpha)}x^{\alpha} \cdot \varphi(\alpha, 2\alpha,x),
\]
after having transformed which series by means of formula (6.), it is
\[
~~~~~~~~~~~~~~~~~~~~~~~~~\frac{x^{\alpha}}{\Pi(\alpha-1)} \int_0^{\infty} (u+u^2)^{\alpha-1} \cdot e^{-ux} \cdot \diff{d}u
\] 
\[
~~~~~~~=\frac{\Pi(2\alpha-2)}{\Pi(\alpha-1)}x^{1-\alpha} \cdot e^{\frac{x}{2}}\cdot \psi(\frac{3}{2}-\alpha,\frac{x^2}{16})+ \frac{\Pi(-2\alpha)}{\Pi(-\alpha)}x^{\alpha} \cdot \psi(\frac{1}{2}+\alpha,\frac{x^2}{16}),
\]
if further $x$ is changed into $4\sqrt{x}$ and $\alpha$ into $\alpha+\frac{1}{2}$, we have by a few reductions
\[
~~~~~~~~17. ~~~~~ \frac{2^{2\alpha+1} \cdot \sqrt{\pi} \cdot x^{\alpha} \cdot e^{-2\sqrt{x}}}{\Pi(\alpha-\frac{1}{2})} \int_0^{\infty} (u+u^2)^{\alpha-\frac{1}{2}} e^{-4u\sqrt{x}} \diff{d}u
\]
\[
~~~~~~~~~~~=\Pi(\alpha-1)\psi(1-\alpha,x)+\Pi(-\alpha-1)x^{\alpha}\psi(1+\alpha,x),
\]
from there it follows by comparision to formula (12.)
\[
\int_0^{\infty} u^{\alpha-1} \cdot e^{-u} \cdot e^{-\frac{x}{u}} \cdot \diff{d}u=\frac{2^{2\alpha+1} \cdot \sqrt{\pi} \cdot x^{\alpha} \cdot e^{-2\sqrt{x}}}{\Pi(\alpha-\frac{1}{2})} \int_0^{\infty} (u+u^2)^{\alpha-\frac{1}{2}} e^{-4u\sqrt{x}} \diff{d}u,
\]
from this formula, or if you like it better, from formula (12.), after having put $\alpha=\frac{1}{2}$, this very simple value of the integral is easily deduced
\[
~~~~~~~~~~~~18.~~~~ \int_0^{\infty} e^{-u^2} \cdot e^{-\frac{x}{u^2}} \cdot \diff{d}u = \frac{\sqrt{\pi}}{2} \cdot e^{-2\sqrt{x}}.
\]
The integrals, the we just found, have muliple applications in Analyisis, for the sake of an example in the integration of the Riccati equation, that, by means of easy substitutions can be changed into the form of equation (9.); I will not sprend more time on these things, but will rather also settle the question about other similar integrals, whose first I chose to be this one:
\[
~~~~~~~~~~~~19.~~~ z= \int_0^{\frac{\pi}{2}} \cos^{\alpha-1} v\cdot \cos(\frac{1}{2}x \tan{v}+\beta v)\diff{d}v.
\] 
I always suppose the quantity $x$ to be positive, because its negative sign can be transferred to the quantity $\beta$. By dfferentiation of the quantity
\[
~~~~~~~~~~~~~~~~~~~~~~~\cos^{\alpha-1} v\cdot \sin(\frac{1}{2}x \tan{v}+\beta v)
\]
it is 
\[
 \diff{d} (\cos^{\alpha-1} v\cdot \sin(\frac{1}{2}x \tan{v}+\beta v))= \cos^{\alpha-2} v\cdot \sin v \cdot \sin(\frac{1}{2}x \tan{v}+\beta v)\diff{d}v
\]
\[
~~~~~~~~~~~~~~~+(\frac{x}{2\cos^2 v}+\beta) \cos^{\alpha-1} v \cos(\frac{1}{2}x \tan{v}+\beta v)\diff{d}v,
\]
and by integrating between the boundaries $v=0$ and $v=\frac{\pi}{2}$
\begin{alignat*}{9}
	&20.~~~~&& 0 && = -(\alpha-1)\int_0^{\frac{\pi}{2}}\cos^{\alpha-2} v\cdot \sin v \cdot \sin(\frac{1}{2}x \tan{v}+\beta v)\diff{d}v&&  &&   &&  &&  &&  \\
	& &&  && +\frac{x}{2}\int_0^{\frac{\pi}{2}}\cos^{\alpha-3} v\cdot \cos(\frac{1}{2}x \tan{v}+\beta v)\diff{d}v &&   &&  &&  &&  \\	
	& &&  && +\beta\int_0^{\frac{\pi}{2}}\cos^{\alpha-1} v\cdot \cos(\frac{1}{2}x \tan{v}+\beta v)\diff{d}v, &&   &&  &&  &&  \\
\end{alignat*}
it is further
\begin{alignat*}{9}
& \frac{\diff{d}z}{\diff{d}x}&&~~~~= \frac{1}{2} \int_0^{\frac{\pi}{2}}\cos^{\alpha-2} v\cdot \sin v \cdot \sin(\frac{1}{2}x \tan{v}+\beta v)\diff{d}v,&&\\
& \frac{\diff{d}^2z}{\diff{d}x^2}&&~~~~= -\frac{1}{4} \int_0^{\frac{\pi}{2}}\cos^{\alpha-3} v\cdot \sin v \cdot \cos(\frac{1}{2}x \tan{v}+\beta v)\diff{d}v,&&\\
\end{alignat*}
therefore
\[
~~~~~~~~~~~z-\frac{\diff{d}^2z}{\diff{d}x^2}= \int_0^{\frac{\pi}{2}}\cos^{\alpha-3} \cdot \cos(\frac{1}{2}x \tan{v}+\beta v)\diff{d}v,
\]
after having substituted which values, equation (20.) is converted into this one
\[
~~~~~~~~~~~21. ~~~~ 0= (x+2 \beta)z+4(\alpha-1)\frac{\diff{d}z}{\diff{d}x}-4x \frac{\diff{d}^2z}{\diff{d}x^2},
\]
whose complete integral is:
\[
~~~~~~~~y= A \psi(\frac{\beta-\alpha+1}{2}, 1-\alpha, x)+Bx^{\alpha} \psi(\frac{\beta+\alpha+1}{2}, 1+\alpha, x),
\]
and because $z= e^{-\frac{x}{2}} \cdot y$, it is
\[
~~~~~~~~~~~22. ~~~~~ \int_0^{\frac{\pi}{2}}\cos^{\alpha-1} \cdot \cos(\frac{1}{2}x \tan{v}+\beta v)\diff{d}v
\]
\[
~~~~= A \cdot \psi(\frac{\beta-\alpha+1}{2}, 1-\alpha, x)+Bx^{\alpha} \psi(\frac{\beta+\alpha+1}{2}, 1+\alpha, x).
\]
The determination of the constant $A$ is easily obtained by putting $x= \infty$, if $a$ is a positive quantity, the determination of the other constant on the other hand requires peculiar artifices; we will obtain both constants by this method: Let us multiply equation (22.) by $x^{\lambda-1} e^{-\frac{x}{2}} \diff{d}x$ and integrate between the boundaries $x=0$ and $x= \infty$, whereafter it is
\begin{alignat*}{9}
	&&&  && \int_0^{\infty} \int_0^{\frac{\pi}{2}}\cos^{\alpha-1} v\cdot x^{\lambda-1} e^{-\frac{x}{2}}  \cos(\frac{1}{2}x \tan{v}+\beta v)\diff{d}v \diff{d}x&&  &&   &&  &&  &&  \\
	& &&  && =A \int_0^{\infty} x^{\lambda-1} \cdot e^{-x} \psi(\frac{\beta-\alpha+1}{2}, 1-\alpha,x)\diff{d}x &&   &&  &&  &&  \\	
	& &&  && =B \int_0^{\infty} x^{\lambda+\alpha-1} \cdot e^{-x} \psi(\frac{\beta+\alpha+1}{2}, 1+\alpha,x)\diff{d}x. &&   &&  &&  &&  \\
\end{alignat*}
The values of all of these integrals can be expressed by known functions, for it is
\[
~~~~~~~~~\int_0^{\infty} x^{c-1} \cdot e^{-x} \cdot \psi(a,b,x)\diff{d}x = \Pi(c-1)F(c,a,b,1),
\]
where $F$ denotes the known hypergeometric series, by which, after having expressed it by the the function $\Pi$, it is
\[
\int_0^{\infty} x^{c-1} \cdot e^{-x} \cdot \psi(a,b,x)\diff{d}x = \frac{\Pi(c-1)\Pi(b-1)\Pi(b-a-c-1)}{\Pi(b-a-1)\Pi(b-c-1)},
\]
further it is
\[
\int_0^{\infty} x^{\lambda-1} \cdot e^{-\frac{x}{2}} \cos{(\frac{1}{2}x \tan{v}+\beta v)}\diff{d}x= 2^{\lambda} \Pi(\lambda-1) \cos^{\lambda} v \cdot \cos{(\lambda+\beta)v},
\]
whose value is expressed by means of the function $\Pi$ in this way
\[
~~~~~~~~~~~~~~~~~~~~~~~~\frac{\pi \cdot \Pi(\lambda-1)\Pi(\alpha+\lambda-1)}{2^{\alpha} \Pi(\frac{\alpha-\beta-1}{2})\Pi(\frac{\alpha+\beta-1}{2}+\lambda)},
\]
after having substituted all which values, equation (23.) is converted into this one:
\[
~~~~~~~~~~~~~~~~~~~~~~~~\frac{\pi \cdot \Pi(\lambda-1)\Pi(\alpha+\lambda-1)}{2^{\alpha} \Pi(\frac{\alpha-\beta-1}{2})\Pi(\frac{\alpha+\beta-1}{2}+\lambda)}
\]
\[
=A \frac{\Pi(\lambda-1)\Pi(-\alpha)\Pi(-\frac{\alpha+\beta+1}{2}-\lambda)}{\Pi(-\frac{\alpha+\beta+1}{2})\Pi(-\alpha-\lambda)}+B\frac{\Pi(\alpha+\lambda-1)\Pi(\alpha)\Pi(-\frac{\alpha+\beta+1}{2}-\lambda)}{\Pi(\frac{\alpha-\beta-1}{2})\Pi(-\lambda)},
\]
this equation is easily reduced to this mre convenient form
\[
~~~~~~~~~\frac{\pi \cdot \cos(\frac{\alpha+\beta}{2})\pi}{2^{\alpha} \Pi(\frac{\alpha-\beta-1}{2})}= \frac{A \cdot \Pi(-\alpha)\sin(\alpha+\lambda)\pi}{\Pi(-\frac{\alpha+\beta+1}{2})}+\frac{B \cdot \Pi(\alpha)\sin \lambda \pi}{\Pi(-\frac{\alpha-\beta-1}{2})},
\]
which, because it holds for any arbitrary value of the quantity $\lambda$, is converted into these two 
\begin{alignat*}{9} 
  && \frac{\pi \cos{\frac{\alpha+\beta}{2}\pi}}{2^{\alpha} \Pi(\frac{\alpha-\beta-1}{2})}&&=& \frac{A \cdot \sin {\alpha \pi} \Pi(-\alpha)}{\Pi(-\frac{\alpha+\beta+1}{2})},&&\\
  &&  -\frac{\pi \sin{\frac{\alpha+\beta}{2}\pi}}{2^{\alpha} \Pi(\frac{\alpha-\beta-1}{2})}&&=& \frac{A \cdot \cos {\alpha \pi} \Pi(-\alpha)}{\Pi(-\frac{\alpha+\beta+1}{2})}&&+\frac{B \cdot \Pi(\alpha)}{\Pi(\frac{\alpha-\beta-1}{2})},\\
\end{alignat*}
from which the values of the constants $A$ and $B$ are easily found
\[
~~~~~A= \frac{\pi \cdot \Pi(\alpha-1)}{2^{\alpha} \Pi(\frac{\alpha-\beta-1}{2})\Pi(\frac{\alpha+\beta-1}{2})}, ~~~~~~~ B= -\frac{\pi \cdot \cos(\frac{\alpha-\beta}{2})\pi}{2^{\alpha} \cdot \sin{\alpha \pi} \Pi(\alpha)},
\]
after having finally substituted which values of the constants in equation (22.), it is
\[
~~~~~~~~~~24. ~~~~ \int_0^{\frac{\pi}{2}} \cos^{\alpha-1} v \cdot \cos(\frac{1}{2}x \tan{v} +\beta v)\diff{d}v
\]
\[
=\frac{\pi \cdot \Pi(\alpha-1)e^{-\frac{x}{2}} \cdot \psi(\frac{\beta-\alpha+1}{2},1-\alpha,x)}{2^{\alpha} \cdot \Pi(\frac{\alpha-\beta-1}{2})\Pi(\frac{\alpha+\beta-1}{2})}-\frac{\pi \cdot \cos{\frac{\alpha-\beta}{2}\pi} \cdot x^{\alpha} \cdot e^{-\frac{x}{2}} \psi(\frac{\beta+\alpha+1}{2},1+\alpha,x)}{2^{\alpha} \sin{\alpha \pi} \Pi(\alpha)}.
\]
Very simple special cases of this formula are:
\begin{alignat*}{9}
&&27. ~~~~ &&\int_0^{\frac{\pi}{2}} \cos^{\alpha-1} v \cdot \cos(x \tan{v} -(\alpha+1)v)\diff{d}v &=\frac{\pi \cdot x^{\alpha} \cdot e^{-x}}{\Pi(\alpha)},&&\\  
&&28. ~~~~ &&\int_0^{\frac{\pi}{2}} \cos^{\alpha-1} v \cdot \cos(x \tan{v} +(\alpha+1)v)\diff{d}v &=0,  \\
\end{alignat*}
of which the one is obtained, after having put $\beta=-\alpha-1$, the other, after having put $\beta=\alpha+1$. From the connected formulas (25.) and (26.) also these ones follow:
\begin{alignat*}{9}
&&27. ~~~~ &&\int_0^{\frac{\pi}{2}} \cos^{\alpha-1} v \cdot \cos(x \tan{v})\cos(\alpha+1)v\cdot \diff{d}v &=\frac{\pi \cdot x^{\alpha} \cdot e^{-x}}{2\Pi(\alpha)},&&\\  
&&28. ~~~~ &&\int_0^{\frac{\pi}{2}} \cos^{\alpha-1} v \cdot \sin(x \tan{v})\sin(\alpha+1)v\cdot \diff{d}v &=\frac{\pi \cdot x^{\alpha} \cdot e^{-x}}{2\Pi(\alpha)}.  \\
\end{alignat*}\\

The formulas (25.) and (26.) agree with the formula, found by the Ill. Laplace, which others later demonstrated in other ways, confer this journal's volume $XIII$, p. 231, where the Cl. Liouville \cite{3}, by the method of differentiation, found for arbitrary parameters
\[
~~~~~~~~~~~~~~~~~~~~~~\int_{-\infty}^{\infty} \frac{e^{\alpha\sqrt{-1}} \cdot \diff{d}\alpha}{(x+\alpha \sqrt{-1})^{\mu}}= \frac{2\pi \cdot e^{-x}}{\Gamma(\mu)}.
\]
Formula (24.) gives another very simple integral, after having put $\beta=\alpha-1$
\[
~~~~~~~~29. ~~~ \int_0^{\frac{\pi}{2}} \cos^{\alpha-1} v \cdot \cos{(x\tan{v}+(\alpha-1)v)}\diff{d}v= \frac{\pi e^{-x}}{2}.
\]
The two series, that are contained in the one part of equation (24.), after having put $\beta$, become $\varphi(\frac{1-\alpha}{2}, 1-\alpha,x)$ and $\varphi(\frac{1+\alpha}{2}, 1+\alpha,x)$, and they can be, by means  of formula (6.), transformed into series of the kind $\psi$. After having done the transformations, if one changes $\alpha$ into $2\alpha$ and $x$ into $4\sqrt{x}$, this formula emerges
\[
~~~~~~~~30. ~~~~~ \frac{2\Pi(\alpha-\frac{1}{2})}{\sqrt{\pi}} \int_0^{\frac{\pi}{2}} \cos^{2\alpha-1} v \cdot \cos(2\sqrt{x} \tan{v})\diff{d}v
\]
\[
~~~~~~~= \Pi(\alpha-1)\psi(1-\alpha,x)+\Pi(-\alpha-1) \cdot x^{\alpha} \cdot \psi(1+\alpha,x),
\]
from this, by comparison to formula (12.), it is
\[
31. ~~~~ \frac{2\Pi(\alpha-\frac{1}{2})}{\sqrt{\pi}} \int_0^{\frac{\pi}{2}} \cos^{2\alpha-1} v \cdot \cos(2\sqrt{x} \tan{v})\diff{d}v= \int_0^{\infty} u^{\alpha-1} \cdot e^{-u} \cdot e^{-\frac{x}{u}} \cdot \diff{d}u.
\]\\

Likewise the connection of the two integrals can be demonstrated, that are contained in the equations (13.) and (24.); for this formula (24.), if $\alpha-\beta$ is put in the place of $\alpha$ and $\alpha+\beta-1$ in the place of $\beta$ and multiplied by $\frac{1}{\pi}\Pi(-\beta) \cdot 2^{\alpha} \cdot e^{\frac{x}{2}} \cdot x^{\beta}$, obtains this form
\[
32. ~~~\frac{2\Pi(-\beta)\cdot e^{\frac{x}{2}}\cdot x^{\beta}}{\pi}\int_0^{\frac{\pi}{2}} (2\cos v)^{\alpha-\beta-1} \cdot \cos(\frac{1}{2}x \tan{v} +(\alpha+\beta-1)v)\diff{d}v
\] 
\[
=\frac{\Pi(\alpha-\beta-1)}{\Pi(\alpha-1)}x^{\beta} \psi(\beta, \beta-\alpha+1,x)+\frac{\Pi(\beta-\alpha-1)}{\Pi(\beta-1)}x^{\alpha} \psi(\alpha, \alpha-\beta+1,x),
\]
after having compared which to formula (13.), it is seen to be
\[
~~~~~~~~~~~~~~~~~~~~33. ~~~~ \int_0^{\infty} \frac{u^{\beta} \cdot e^{-ux} \cdot \diff{d}u}{(1+u)^{\alpha}}
\]
\[
~~~~~~=\frac{2 \cdot e^{x}{2}}{\sin{\beta \pi}}\int_0^{\frac{\pi}{2}} (2\cos v)^{\alpha-\beta-1} \cdot \cos(\frac{1}{2}x \tan{v} +(\alpha+\beta-1)v)\diff{d}v,
\]
moreover, since the one part of equation (32.) can be transformed by means of formula (7.), it is
\[
34. ~~~~ \frac{2\Pi(-\beta)e^{\frac{x}{2}}\cdot x^{\beta}}{\pi}=\int_0^{\frac{\pi}{2}} (2\cos v)^{\alpha-\beta-1} \cdot \cos(\frac{1}{2}x \tan{v} +(\alpha+\beta-1)v)\diff{d}v= \chi(\alpha, \beta,x).
\]
We will also treat this more general integral in the same way
\[
~~~~~~~~~~~~~~y= \int_0^{\frac{\pi}{2}} \sin^{\alpha-1} v \cdot \cos^{\beta-1} \cdot \cos(x\tan{v}+\gamma v)\diff{d}v
\]
and we will chose the cases, in which it can be expressed by the series mentioned above. We also suppose the quantity $x$ always to be positive in this integral, because it is possible, to transfer its negative sign to the quantity $\gamma$. By differentiating the formula $\sin^{\alpha} v \cdot \cos^{\beta} \cdot \cos(x\tan{v}+\gamma v)$, and integrating from $u=0$ to $u=\frac{\pi}{2}$ thereafter, it is
\begin{alignat*}{9}
	& && 0 && = \alpha \int_0^{\frac{\pi}{2}} \sin^{\alpha-1} v \cdot \cos^{\beta+1} \cdot \cos(x\tan{v}+\gamma v)\diff{d}v&&  &&   &&  &&  &&  \\
	& &&  && -\beta \int_0^{\frac{\pi}{2}} \sin^{\alpha+1} v \cdot \cos^{\beta-1} \cdot \cos(x\tan{v}+\gamma v)\diff{d}v &&   &&  &&  &&  \\	
	& &&  && -x \int_0^{\frac{\pi}{2}} \sin^{\alpha} v \cdot \cos^{\beta-2} \cdot \sin(x\tan{v}+\gamma v)\diff{d}v &&   &&  &&  &&  \\
	& &&  && -\gamma \int_0^{\frac{\pi}{2}} \sin^{\alpha} v \cdot \cos^{\beta} \cdot \sin(x\tan{v}+\gamma v)\diff{d}v, &&   &&  &&  &&  \\
\end{alignat*}
from this equation, if the integrals are expressed by $y$ and its differentials, this differential equation of third order is easily deduced
\[
~~~~~~~~~~ 35. ~~~ 0= \alpha y +(\gamma+x)\frac{\diff{d}y}{\diff{d}x}+(\beta-2)\frac{\diff{d}^2y}{\diff{d}x^2}-x\frac{\diff{d}^3y}{\diff{d}x^3},
\]
if one puts 
\[
~~~~~~~~~~~~~~~~36. ~~~~ A_0+A_1x+A_2x^2+A_3x^3+\cdots, 
\]
the conditional equations are easily found, that have to hold between the coefficients of this series, that this series satisfies the differential equation
\begin{alignat*}{9}
 & ~~~~~~~&& \alpha &&A_0 &&+\gamma \cdot 1 \cdot A_1 && -1 \cdot 2 \cdot (2-\beta)A_2, &&\\
 & && (\alpha +1)&& A_1 &&+\gamma \cdot 2 \cdot A_2 && -2 \cdot 3 \cdot (3-\beta)A_3, &&\\
\end{alignat*}
and in general
\[
~~~ 37. ~~~~~ (\alpha+k)A_k+\gamma \cdot (k+1) A_{k+1}- (k+1)(k+2)(k+2-\beta)A_{k+2}.
\]
If one puts in the same way 
\[
~~~~~~~~~~~~~~ 38. ~~~~ y= x^{\beta}(B_0+B_1x+B_2x^2+B_3x^3+\cdots),
\]
one finds these relations for the coefficients
\begin{alignat*}{9}
 & ~~~~~~~&& &&+\gamma \cdot \beta \cdot B_0 &&-\beta(\beta+1) \cdot 1 \cdot B_1, &&\\
 & && (\alpha +\beta)B_0&& +\gamma(\beta+1)B_1 &&+(\beta+1)(\beta+2) \cdot 2 \cdot B_2, &&\\
\end{alignat*}
and in general 
\[
39. ~~~ (\alpha+\beta+k)B_k+\gamma(\beta+k+1)B_{k+1}-(\beta+k+1)(\beta+k+2)(k+2)B_{k+2},
\]
from this it is clear, that the complete integral of equation (35.) is
\[
~~~~~ 40. ~~~ y= A_0+A_1x+A_2x^2+\cdots+x^{\beta}(B_0+B_1x+B_2x^2+\cdots),
\]
for, by the equations (37.), two of the quantities $A_0$,$A_1$,$A_2$ etc, and by the equations (39.) one of the quantities $B_0$,$B_1$,$B_2$ etc. remain arbitrary, so that this integral contains three arbitrary constants. Therefore, if the integral mentioned above is resubstituted, it is
\[
~~~~~~~~~~~~~~~~ \int_0^{\frac{\pi}{2}} \sin^{\alpha-1} v \cdot \cos^{\beta-1} \cdot \cos(x\tan{v}+\gamma v)\diff{d}v
\]
\[
~~~~~~~~~~~=A_0+A_1x+A_2x^2+\cdots+x^{\beta}(B_0+B_1x+B_2x^2+\cdots).
\]
From the relations of the coefficients it is easily seen, that these series and this general integral are higher transcendentals than those, that we undertake to treat here; but they nevertheless coincide with those in certain special cases. At first, if we suppose $\gamma= \alpha+\beta$, it follows from the equations (39.)
\begin{alignat*}{9}
 &~~~~~~~~~&& B_1 &&=\frac{\alpha+\beta}{1(1+\beta)}B_0, &&\\
 & && B_2 &&=\frac{(\alpha+\beta)(\alpha+\beta+1)}{1\cdot 2(1+\beta)(2+\beta)}B_0, &&\\
 & && B_3 &&=\frac{(\alpha+\beta)(\alpha+\beta+1)(\alpha+\beta+2)}{1\cdot 2 \cdot 3(1+\beta)(2+\beta)(3+\beta)}B_0, &&\\
 & &&\text{etc.}&& ~~~~~~~~~~~~\text{etc.}&&\\
\end{alignat*}
Further, if $\beta$ is positive, after having put $x=0$, it follows from equation (41.)
\[
A_0= \int_0^{\frac{\pi}{2}} \sin^{\alpha-1} v \cdot \cos^{\beta-1} v \cdot \cos((\alpha+\beta)v \cdot \diff{d}v= \frac{\cos{\frac{\alpha \pi}{2}} \Pi(\alpha-1)\Pi(\beta-1)}{\Pi(\alpha-\beta-1)},
\]
if equation (41.) is differentiated with respect to $x$ in the same way and one puts $x=0$ afterwards, it is
\[
A_1= -\int_0^{\frac{\pi}{2}} \sin^{\alpha} v \cdot \cos^{\beta-2} v \cdot \sin((\alpha+\beta)v \cdot \diff{d}v= -\frac{\cos{\frac{\alpha \pi}{2}} \Pi(\alpha)\Pi(\beta-2)}{\Pi(\alpha-\beta-1)},
\]
therefore it is
\[
~~~~~~~~~~~~~~~~~~~~~~~~~~~~~~~A_1=\frac{\alpha}{1(1-\beta)}A_0,
\]
from there it easily follows from the equations (37.)
\begin{alignat*}{9}
 &~~~~~~~~~~~~~~~~~~~~~~~~~&& A_2 &&=\frac{\alpha(\alpha+1)A_0}{1 \cdot 2(1-\beta)(2-\beta)}, &&\\
 & && A_1 &&=\frac{\alpha(\alpha+1)[\alpha+2)A_0}{1 \cdot 2 \cdot 3(1-\beta)(2-\beta)(3-\beta)}, &&\\
 & &&\text{etc.}&& ~~~~~~~~~~~~\text{etc.}&&\\
\end{alignat*}
Therefore these two series, by which we expressed our integral, belong to the class of series, that we denoted by $\varphi$ above, in this case $\gamma=\alpha+\beta$, and formula (41.) is converted into this one:
\[
~~~~~~~~~~~~~~~~ \int_0^{\frac{\pi}{2}} \sin^{\alpha-1} v \cdot \cos^{\beta-1} \cdot \cos(x\tan{v}+(\alpha+\beta) v)\diff{d}v
\]
\[
~~~=\frac{\cos{\frac{\alpha \pi}{2}}\Pi(\alpha-1)\Pi(\beta-1)}{\Pi(\alpha+\beta-1)} \varphi(\alpha, 1-\beta,x)+B_0x^{\beta}\varphi(\alpha+\beta,1+\beta,x).
\]
In the determination of the constant $B_0$ we will use the same method as above in the determination of the constants of equation (22.). By multiplying by $x^{\lambda-1} \cdot e^{-x} \diff{d}x$ and integrating between the boundaries $0$ and $\infty$, it is
\[
~~~~~~~~~~\Pi(\lambda-1) \int_0^{\frac{\pi}{2}} \sin^{\alpha-1} v \cdot \cos^{\beta+\lambda-1} \cdot \cos(\alpha+\beta+\lambda) v\diff{d}v
\]
\[
=\frac{\cos{\frac{\alpha \pi}{2}}\Pi(\alpha-1)\Pi(\beta-1)\Pi(\lambda-1)}{\Pi(\alpha+\beta-1)}F(\lambda, \alpha, 1-\beta,1)+ B_0 \Pi(\beta+\lambda-1)F(\lambda+\beta, \alpha+\beta, 1+\beta,1),
\]
and after having expressed the hypergeometric series along with the integral by the function $\Pi$, it is 
\[
~~~~~~~~~~~~~~~~~~~\frac{\cos{\frac{\alpha \pi}{2}}\Pi(\lambda-1)\Pi(\alpha-1)\Pi(\beta+\lambda-1)}{\Pi(\alpha+\beta+\lambda-1)}
\]
\[
= \frac{\cos{\frac{\alpha \pi}{2}}\Pi(\lambda-1)\Pi(\alpha-1)\Pi(\beta+\lambda-1)\Pi(-\beta)\Pi(-\beta-\alpha-\alpha)}{\Pi(\alpha+\beta+\lambda-1)\Pi(-\alpha-\beta)\Pi(-\beta-\lambda)}
\]
\[
~~~~~~~~~~~~~~~~+B_0\frac{\Pi(\beta+\lambda-1)\Pi(\beta)\Pi(-\beta-\alpha-\lambda)}{\Pi(-\alpha)\Pi(-\beta)},
\]
after some reductions the quantity $\lambda$, because it has to, vanishes completely, and this very simple value of the constant $B_0$ emerges
\[
~~~~~~~~~~~~~~~~~~~~~~~~~~~B_0 = \cos{\frac{\alpha \pi}{2}}\Pi(-\beta-1),
\]
after having finally substituted which value, we have
\[
~~~~~~ 42. ~~~ \int_0^{\frac{\pi}{2}} \sin^{\alpha-1} v \cdot \cos^{\beta-1} v \cdot \cos(x\tan{v}+(\alpha+\beta)v)\diff{d}v
\]
\[
=\frac{\cos{\frac{\alpha \pi}{2}}\Pi(\alpha)\Pi(\beta-1)}{\Pi(\alpha+\beta-1)}\varphi(\alpha,1-\beta,x)+x^{\beta}\cos{\frac{\alpha \pi}{2}}\Pi(-\beta-1)\varphi(\alpha+\beta,1+\beta,x).
\]
and after having compared these formulas to each other, one sees the connection of the two integrals
\begin{alignat*}{9}
& 43. ~~~~ && ~&& \cos{\frac{\alpha \pi}{2}}\int_0^{\frac{\pi}{2}} \sin^{\alpha-1} v \cdot \cos^{\beta-1} v \cdot \sin(x\tan{v}+(\alpha+\beta)v)\diff{d}v\\
& &&=&& \sin{\frac{\alpha \pi}{2}}\int_0^{\frac{\pi}{2}} \sin^{\alpha-1} v \cdot \cos^{\beta-1} v \cdot \cos(x\tan{v}+(\alpha+\beta)v)\diff{d}v,\\
\end{alignat*}
which formula can also be exhibited in this way
\[
~~~ 44. ~~~ \int_0^{\frac{\pi}{2}} \sin^{\alpha-1} v \cdot \cos^{\beta-1} v \cdot \sin(x\tan{v}+(\alpha+\beta)v-\frac{\alpha \pi}{2})\diff{d}v=0
\]
The special case of formula (42.), in which $\alpha=0$, is worth to be noted
\[
~~~~~~~~~~~~~~ 45. ~~~~ \int_0^{\frac{\pi}{2}} \frac{\cos^{\beta-1} \cdot \sin(x \tan{v}+\beta v}{\sin{v}}\diff{d}v= \frac{\pi}{2},
\]
of which the Cl. Liouville found the special case, corresponding to the value $x=0$, in this journal, volume $XIII$, page 232 \cite{3}. Moreover, having compared the formulas (42.) and (13.), one sees the connection of this integral to those, that we treated above, without any difficulty
\[
~~~~~~~~~46. ~~~~~\frac{\cos{\frac{\alpha \pi}{2}} \Pi(\alpha-1)}{\Pi(\alpha+\beta-1)}x^{\beta} \int_0^{\infty} \frac{u^{\alpha+\beta-1}\cdot e^{-ux}\diff{d}u}{(1+u)^{\alpha}}
\]
\[
~~~~~=\int_0^{\frac{\pi}{2}} \sin^{\alpha-1} v \cdot \cos^{\beta-1} v \cdot \cos(x\tan{v}+(\alpha+\beta)v)\diff{d}v.
\]
Another case, in which the series of formula (41.) are converted into series denoted by the character $\psi$, is $\gamma= -\alpha-\beta$, for in this case formula (41.) is easily found by the same method as above, to be converted into this one:
\[
~~~~~~~~~=\int_0^{\frac{\pi}{2}} \sin^{\alpha-1} v \cdot \cos^{\beta-1} v \cdot \cos(x\tan{v}-(\alpha+\beta)v)\diff{d}v
\]
\[
= \frac{\cos{\frac{\alpha \pi}{2}}\Pi(\alpha-1) \Pi(\beta-1)}{\Pi(\alpha+\beta-1)}\varphi(\alpha, 1-\beta, -x)+B_0x^{\beta} \varphi(\alpha+\beta, 1+\beta, -x),
\]
but in this case the constant $B_0$ obtains another value, that we find by muliplying by $x^{\alpha+\beta} \cdot e^{-x} \diff{d}x$ and by integrating between the boundaries $x=0$ and $x=\infty$, after having done those integrations, it is
\[
~~~~~~~~~~~\Pi(\alpha+\beta-1) \int_0^{\frac{\pi}{2}} \sin^{\alpha-1} v \cdot \cos^{\alpha+2\beta-1}v \cdot \diff{d} v
\]
\[
~~~~~~~=\cos{\frac{\alpha \pi}{2}} \Pi(\alpha-2) \Pi(\beta-1)F(\alpha+\beta, \alpha, 1-\beta, -1)
\]
\[
~~~~~~~~~~+B_0 \Pi(\alpha+2\beta-1)F(\alpha+2\beta, \alpha+\beta, 1+\beta, -1),
\]
these hypergeoemtric series, whose forth element is $=-1$, can also be expressed by the function $\Pi$ by means of the formula
\[
~~~~~~~~~~~~~~~F(\alpha, \beta, \alpha-\beta+1, -1) = \frac{2^{-\alpha} \sqrt{\pi} \Pi(\alpha-\beta)}{\Pi(\frac{\alpha}{2}-\beta)\Pi(\frac{\alpha-1}{2})},
\]
that I proved in the commentary about the hypergeometric series in this journal volume $XV$ page 135 \cite{2}. From this, if that integral and the hypergeometric series are expressed by the function $\Pi$, it emerges after certain easy reductions:
\[
~~~~~~~~~~~~~~~~~~~~~~~B_0= \cos(\frac{\alpha}{2}+\beta)\pi \cdot \Pi(-\beta-1),
\]
and after having substituted the value of the constant, it is:
\[
~~~ 47. ~~~~\int_0^{\frac{\pi}{2}} \sin^{\alpha-1} v \cdot \cos^{\beta-1} v \cdot \cos(x\tan{v}-(\alpha+\beta)v)\diff{d}v
\]
\[
= \frac{\cos{\frac{\alpha \pi}{2}}\Pi(\alpha-1) \Pi(\beta-1)}{\Pi(\alpha+\beta-1)}\varphi(\alpha, 1-\beta, -x)+x^{\beta}\cos(\frac{\alpha}{2}+\beta)\pi \cdot \Pi(-\beta-1)  \varphi(\alpha+\beta, 1+\beta, -x).
\]
A similar formula is deduced from this one, by changing $\alpha$ into $\alpha-1$, $\beta$ into $\beta+1$ and differentiating with respect to the variable $x$
\[
~~~ 48. ~~~~\int_0^{\frac{\pi}{2}} \sin^{\alpha-1} v \cdot \cos^{\beta-1} v \cdot \sin(x\tan{v}-(\alpha+\beta)v)\diff{d}v
\]
\[
= -\frac{\sin{\frac{\alpha \pi}{2}}\Pi(\alpha-1) \Pi(\beta-1)}{\Pi(\alpha+\beta-1)}\varphi(\alpha, 1-\beta, -x)-x^{\beta}\sin(\frac{\alpha}{2}+\beta)\pi \cdot \Pi(-\beta-1)  \varphi(\alpha+\beta, 1+\beta, -x).
\]
These formulas (47.) and (48.) can easily be combined in two ways like this, that they obtain these simpler forms:
\[
~~~ 49. ~~~~\int_0^{\frac{\pi}{2}} \sin^{\alpha-1} v \cdot \cos^{\beta-1} v \cdot \sin(x\tan{v}-(\alpha+\beta)v+(\frac{\alpha}{2}+\beta)\pi)\diff{d}v
\]
\[
~~~~~~~~~~~~~~~~~~~~~~= \frac{\pi \Pi(\alpha-1)\varphi(\alpha, 1-\beta,-x)}{\Pi(-\beta)\Pi(\alpha+\beta-1)},
\]
\[
~~~ 50. ~~~~\int_0^{\frac{\pi}{2}} \sin^{\alpha-1} v \cdot \cos^{\beta-1} v \cdot \sin(x\tan{v}-(\alpha+\beta)v+\frac{\alpha \pi}{2})\diff{d}v
\]
\[
~~~~~~~~~~~~~~~~~~~~~~= \frac{\pi x^{\beta}}{\Pi(\beta)}\varphi(\alpha+\beta, 1+\beta, -x).
\]
In all these integrals, that were treated here, as we already mentioned earlier, $x$ has to be a positive quantity, but if $x$ is supposed to be negative, all sums found would be wrong; in this the integral of equation (50.) is worth to be noted, that, for positive $x$, is equal to that series, but vanishes for negative $x$, confer equation (44.).\\[5mm]

Legnica, in the month of April, 1837

\end{document}